\documentclass{amsart}
\usepackage[latin1]{inputenc}
\usepackage{amsmath}
\usepackage{amsthm}
\usepackage{amsfonts}
\usepackage{cases}
\usepackage{latexsym}
\usepackage{stmaryrd}
\usepackage{amsfonts}
\usepackage{amsmath,amssymb,amscd,bbm,amsthm,mathrsfs,dsfont}
\usepackage[T1]{fontenc}
\usepackage{float}
\usepackage[all]{xy}
\usepackage{tikz}
\usepackage{calrsfs}
\usepackage{bbm}
\usepackage{color}
\usepackage[notref,notcite]{showkeys}
\usepackage{stmaryrd}

\usepackage{fancyhdr}
\usepackage{amsxtra,ifthen}
\usepackage{verbatim}

\usepackage[numbers,sort&compress]{natbib}
\usepackage[pagebackref]{hyperref}
\usepackage{hypernat}

\newtheorem{theorem}{Theorem}[section]
\newtheorem{lemma}[theorem]{Lemma}
\newtheorem{question}[theorem]{Question}
\newtheorem{e-proposition}[theorem]{Proposition}
\newtheorem{corollary}[theorem]{Corollary}
\newtheorem{e-definition}[theorem]{Definition\rm}

 \DeclareMathOperator{\gld}{gl.dim}
\DeclareMathOperator{\Rat}{Rat}

 \DeclareMathOperator{\End}{End}
 
\DeclareMathOperator{\Hom}{Hom}

\def \lim{\hbox{\lower0.9ex\hbox{$\buildrel{\lower1.1ex\hbox{lim}}
\over{\textstyle \leftarrow}$}}\,}

\begin{document}

\title[The Iwasawa Algebra $\Omega_G$ and Its Dual Artin Coalgebra]
{The Iwasawa Algebra $\Omega_G$ and Its Dual Artin Coalgebra}

\author{Zheng Fang and Feng Wei}

\address{Fang: School of Mathematics and Statistics, Beijing Institute of
Technology, Beijing, 100081, P. R. China}

\email{514139626@qq.com}

\address{Wei: School of Mathematics and Statistics, Beijing Institute of
Technology, Beijing, 100081, P. R. China}

\email{daoshuo@hotmail.com}\email{daoshuowei@gmail.com}

\begin{abstract}
For any compact $p$-adic Lie group $G$,
the Iwasawa algebra $\Omega_G$ over finite field $\mathbb{F}_p$
is a complete noetherian semilocal algebra. It is shown that $\Omega_G$
is the dual algebra of an artinian coalgebra $C$. We induce
a duality between the derived category $\mathcal{D}^b_{fg}(_{\Omega_G}\mathcal{M})$
of bounded complexes of left $\Omega_G$-modules with finitely generated
cohomology modules
and the derived category $\mathcal{D}^b_{qf}(^C\mathcal{M})$
of bounded complexes of left $C$-comodules with quasi-finite
cohomology comodules.

\vskip 0.5\baselineskip

\noindent{\bf L'alg\`ebre d'Iwasawa $\Omega_G$ et Son Coalg\`ebre
d'Artin Duale.}

\noindent{\bf R\'esum\'e}

\noindent Pour tout groupe de Lie $G$ $p$-adique compact, l'alg\`ebre
d'Iwasawa $\Omega_G$ sur un corps fini $\mathbb{F}_p$ est
une alg\`ebre noetherienne, semilocale et compl\`ete. Il est montr\'e que $\Omega_{G}$ est l'alg\`ebre duale d'une coalg\`ebre artinienne $C$. On induit une dualit\'e entre la cat\'egorie  d\'eriv\'ee
$\mathcal{D}^b_{fg}(_{\Omega_G}\mathcal{M})$ des complexes born\'es
de $\Omega_G$-modules \`a gauche avec des modules cohomologie de type fini
et la cat\'egorie  d\'eriv\'ee
$\mathcal{D}^b_{qf}(^C\mathcal{M})$ des complexes born\'es de
$C$-comodules \`a gauche dont leur comodules de cohomologie sont quasi-finis. 

\end{abstract}



%\thanks{????}

\date{02 October 2016}

\maketitle

%\tableofcontents

\section*{Version fran\c{c}aise abr\'eg\'ee}
Dans \cite{Yekutieli}, Yekutieli introduit le concept de 
complexe dualisant non
commutatif. Avec les complexes dualisants, on peut trouver 
plusieurs propri\'et\'es assez remarquables sur les alg\`ebres 
noeth\'erien semilocale et compl\`ete (voir \cite{Wei}, \cite{WuZhang1}
and \cite{WuZhang2}). Pour tout groupe de Lie $G$ $p$-adique
compact, l'alg\`ebre d'Iwasawa $\Omega_G$ sur le corps fini
$\mathbb{F}_p$ est semilocale noeth\'erien  compl\`ete (voir 
\cite{ArdakovBrown}, \cite{Lazard}). Dans \cite{Wei}, on a utilis\'e 
le complexe dualisant pour   \'etudier les propri\'et\'es homologiques 
de l'alg\`ebre d'Iwasawa $\Omega_G$. Il a \'et\'e montr\'e que la formule
d'Auslander-Buchsbaum, le th\'eor\`eme de Bass et le th\'eor\`eme
des $\ll$non trous$\gg$ maintenez pour l'alg\`ebre d'Iwasawa
$\Omega_G $. Une version duale correspondant de ces r\'esultats a
\'et\'e obtenue par la th\'eorie de dualit\'e de Morita. De plus,  
il a \'et\'e prouv\'e que l'alg\`ebre d'Iwasawa $\Omega_G$ est
auto-duale au sens de Morita. On renvoit le lecteur \`a \cite{Wei} pour cettes r\'esultats et pour plus de d\'etails.

D'autre part, il n'est pas difficile de montrer que le radical de Jacobson $J(\Omega_{G})$ de l'alg\`ebre d'Iwasawa $\Omega_G$ est cofini, c'est \`a dire,
l'alg\`ebra quotient $\Omega_G/J(\Omega_G)$ de $\Omega_G$ avec
l'\'egard de ses Jacobson radical $J(\Omega_G)$ est de dimension
finie. Heyneman
et Radford (\cite{HeynemanRadford}) ont d\'emontr\'e qu'une 
alg\`ebre noeth\'erien    
compl\`ete dont le radical de Jacobson est cofini est l'alg\`ebre 
duale d'une coalg\`ebre artinienne. Cela nous motive \`a consid\'erer
quelques propri\'et\'es de dualit\'e entre la cat\'egorie de $\Omega_G$-modules et la cat\'egorie de $C$-modules. Dans cet article, on induit une dualit\'e entre la cat\'egorie deriv\'ee $\mathcal{D}^b_{fg}(_{\Omega_G}\mathcal{M})$ des complexes born\'es de $\Omega_G$-modules \`a gauche \`a cohomologie de type fini, et la cat\'egorie d\'eriv\'ee $\mathcal{D}^b_{qf}(^C\mathcal{M})$ des
complexes born\'es de $C$-comodules \`a gauche dont leurs comodules de cohomologie sont quasi-finis.

\setcounter{section}{-1}
\section{Introduction}
\label{xxsec0}

In the recent years, there has been an increasing interest in the
noncommutative Iwasawa algebras of compact $p$-adic Lie groups. The
main motivation to study the Iwasawa algebras $\Lambda_G$ and
$\Omega_G$ of a compact $p$-adic Lie group $G$
is due to their connections with number theory and arithmetic
algebriac geometry (see \cite{Coates},
\cite{CoatesSchneiderSujatha}, \cite{CoatesFukayaKatoSujathaVenjakob},
\cite{OchiVenjakob},
\cite{Venjakob3}, \cite{Venjakob4}, \cite{Venjakob5}). Throughout, let $p$ be a fixed
prime integer, $\mathbb{Z}_p$ be the ring of $p$-adic integers and
$\mathbb{F}_p$ be the field of $p$ elements.
Given a compact $p$-adic Lie group, the \textit{Iwasawa algebra}
of $G$ over $\mathbb{Z}_p$ is the completed group algebra
$$
\Lambda_G:=\underleftarrow{{\rm lim}}\mathbb{Z}_p[G/N],
$$
where the inverse limit is taken over the open normal subgroups $N$
of $G$. Closely related to $\Lambda_G$ is the \textit{Iwasawa algebra}
of $G$ over $\mathbb{Z}_p$, which is defined as
$$
\Omega_G:=\underleftarrow{{\rm lim}}\mathbb{F}_p[G/N].
$$
It is well-known that $\Lambda_G$ and
$\Omega_G$ are complete
noetherian semilocal algebras and are in general noncommutative.
These algebras associated with certain topological setting
were defined and studied
by Lazard in his seminal paper \cite{Lazard} at first. We would like to
point out that the Iwasawa algebras $\Lambda_G$ and $\Omega_G$, despite its
great interest in number theory and arithmetic, seems to have
been neglected as two concrete examples for the application
of noncommutative methods. Some papers are specially contributed to
the ring-theoretic (or module-theoretic) and homological
properties of Iwasawa algebras (see \cite{ArdakovBrown}, \cite{ArdakovWeiZhang1},
\cite{ArdakovWeiZhang2}, \cite{CoatesSchneiderSujatha}, \cite{Venjakob1},
\cite{Venjakob2}, \cite{Venjakob4}, \cite{Wei}).
We refer to \cite{ArdakovBrown}
and \cite{Lazard} for the basic properties of $\Lambda_G$ and
$\Omega_G$.

%The noncommutative dualizing complex, which was introduced by 
%Yekutieli in \cite{Yekutieli}, provides a perspective and powerful tool to
%the study of noncommutative algebras. In order to determine 
%a dualizing complex over a graded algebra, Yekutieli posed the notion
%of balanced dualizing complex. But, the terminology `balanced' make no sense
%for a general non-graded algebra. In this case, we must used the so-called
%rigid dualizing complex, which is due to Van den Bergh \cite{VandenBergh}. 
%Van den Bergh's results were extended to complete noetherian semilocal
%algebras. A number of elegant properties of complete noetherian semilocal
%were discovered via dualizing complexes (see \cite{WuZhang1}, \cite{WuZhang2}). 
%The author of this article 
%applied some of these properties to describe homological behavior of Iwasawa algebras
%\cite{Wei}. He obtained the Auslander-Buchsbaum formula, the Bass's theorem and no-holes %theorem for noetherian modules over $\Omega_G$, and the dual versions for their Aritian %modules. Moreover, it 
%was shown that $\Omega_G[d]$ is an Auslander, pre-balanced, 
%Cdim-symmetric, bifinite dualizing complex over $\Omega_G$, where $d$ is the injective
%dimension of $\Omega_G$. 

Let $G$ be a compact
$p$-adic analytic group and $N$ be a closed normal subgroup of
$G$. Passman in \cite{Passman} showed us 
that if $R = S \ast G$ is a crossed product of a ring $S$ with a finite
group $G$, then $J(S) \ast G$ is contained in $J(R)$. Thus the
augmentation ideal $w_{N,G} =
\ker(\Omega_G \to \mathbb{F}_p[G/N])=J(\Omega_N) \ast G$ is
contained in $J(\Omega_G)$. On the other hand, $\Omega_G /
w_{N,G}\cong \mathbb{F}_p[G/N]$ is finite dimensional. Hence
$\Omega_G/J(\Omega_G)$ is also finite dimensional, being a
quotient of $\Omega_G / w_{N,G}$. In this situation, we say that 
$\Omega_G$ has \textit{cofinite Jacobson radical}.
Heyneman and Radford \cite{HeynemanRadford} observed that any noetherian complete
algebra $A$ with cofinite Jacobson radical is the dual
algebra of an artinian coalgebra $C$. Therefore we conclude that
for any compact $p$-adic analytic group $G$, the Iwasawa algebra
$\Omega_G$ is the dual algebra of an artinian coalgebra $C$ (see Propsotion \ref{xxsec1.1} and Corollary \ref{xxsec1.2}).

The main purpose of this article is to describe certain duality properties 
between the category of $\Omega_G$-modules and the category of $C$-comodules.

\section{The Main Results}
\label{xxsec1}

Let us first recall some basic definitions and existing facts.
A subspace $W$ of a vector space $V$ over a field $K$ is
called \textit{cofinite} in $V$, or simply \textit{cofinite},
when the quotient space $V/W$ is finite dimensional.

Let $A$ be an algebra over a field $K$ and $M$ be an $A$-module. 
$M$ is said to be \textit{almost noetherian}
if every cofinite submodule of $M$ is finitely generated. 
$A$ is called \textit{left almost noetherian}
if it is almost noetherian as a left
$A$-module, i.e., every cofinite left ideal of $A$ is
finitely generated. $A$ is an \textit{almost noetherian} algebra
if $A$ is both left and right almost noetherian. We refer
the reader to the reference \cite{HeynemanRadford} for
basic properties of almost noetherian modules and
those of almost algebras.

Some topological preliminaries are necessary for our
later work. Suppose that $V$ is a vector space over a field $K$.
We will regard the dual space $V^*$
of al1 linear functionals on $V$ as a topological vector space
with the weak-$*$ topology; that is, $V^*$ is given the least fine
topology such that the vectors of $V$ induce
continuous functionals (where we regard $V$ as embedded in $V^{**}$
and give the discrete topology to the scalars $k$).
The closed (respectively, open) linear
subspaces of $V^*$ are then the annihilators $W^\bot$ of arbitrary
(respectively, finite-dimensional) subspaces $W$ of $V$.

Let $A$ be an algebra over a field $K$. If $I$ is a subspace of
$A$, then
$$
I^\bot=\left\{\hspace{2pt}f\in A^*\hspace{2pt}|\hspace{2pt} f(I)=0
\hspace{2pt}\right\}.
$$
Recall \cite[Section 6]{Sweedler} that
$$
A^0=\left\{\hspace{2pt}I^\bot\hspace{2pt}|\hspace{2pt} I \hspace{2pt} {\rm runs}
\hspace{2pt} {\rm over} \hspace{2pt} {\rm the} \hspace{2pt}
{\rm cofinite} \hspace{2pt} {\rm ideals} \hspace{2pt}
{\rm of} \hspace{2pt} A \right\}
$$
inherits a natural coalgebra structure from the transpose of the
algebra structure map $A\bigotimes A \longrightarrow A$. Now let
$C$ be a coalgebra. If $M$ is a finite dimensional
$C^*$-module and $\pi: C^*\longrightarrow {\rm End} \hspace{2pt}M$
is the corresponding algebra map, then $I_M={\rm Ker}\hspace{2pt}\pi$
is a cofinite ideal of $C^*$. By \cite[Lemma 1.2 (e)]{HeynemanRadford} it
follows that $(\pi^*)=I_M^\bot$, the annihilator of $I_M$
in $({C^*})^*$, so that the image $(\pi^*)$ is contained in $(C^*)^0$.

A module $M$ is said to be \textit{rational} if the image $(\pi^*)$
is contained in $C$ (of course $C$ is a subspace of $A^0$). Since we have
seen in \cite[Lemma 1.2 (b)]{HeynemanRadford} that a cofinite subspace $I$
of $A$ is closed when $I^\bot$, the annihilator of $I$ in $A^*$, 
is actually contained in $C$, it follows that $M$ is a rational
$A$-module when $I_M={\rm Ker} \hspace{2pt} \pi$ is closed in $A$.

Let $C_0$ be the sum of the simple subcoalgebras of a coalgebra
$C$. We inductively define an increasing coalgebra filtration by
$$
C_n=C_0\wedge C_{n-1}=\left\{\hspace{2pt}x\in C \hspace{2pt}|\hspace{2pt}
\Delta x\in C_0\otimes C+C\otimes C_{n-1} \hspace{2pt}\right\}.
$$
$C$ is said to be \textit{finite type} when $C_1=C_0\wedge C_0$
is finite dimensional.

\begin{e-proposition}\cite[Proposition 4.3.1]{HeynemanRadford}
\label{xxsec1.1}
Let $A$ be a left almost noetherian algebra with cofinite
Jacobson radical. Then $C=A^0$ is a coalgebra of finite type,
and the algebra map $A\longrightarrow C^*$ can be identified
with the map $A\longrightarrow \widehat{A}$ (completion with respect
to the $J$-adic topology). In particular $C^*=\widehat{A}$
is again an almost noetherian algebra with cofinite Jacobson radical.
\end{e-proposition}

This proposition means that any complete noetherian algebra $A$ 
with cofinite Jacobson
radical is the dual algebra of an artinian coalgebra $C$.
Moreover, the artinian coalgebra $C$ is of finite type.
We will now apply the above results to the Iwasawa algebra
$\Omega_G$, which is our work context. 
%Let $G$ be a compact
%$p$-adic analytic group and $N$ be a closed normal subgroup of
%$G$. Passman in \cite{Passman} has proved
%that if $R = S \ast G$ is a crossed product of a ring $S$ with a finite
%group $G$, then $J(S) \ast G$ is contained in $J(R)$. Thus the
%augmentation ideal $w_{N,G} =
%\ker(\Omega_G \to \mathbb{F}_p[G/N])=J(\Omega_N) \ast G$ is
%contained in $J(\Omega_G)$. On the other hand, $\Omega_G /
%w_{N,G}\cong \mathbb{F}_p[G/N]$ is finite dimensional. Hence
%$\Omega_G/J(\Omega_G)$ is also finite dimensional, being a
%quotient of $\Omega_G / w_{N,G}$. By the
%above discussion we know that the Iwasawa algebra $\Omega_G$
%is an almost noetherian algebra with cofinite Jacobson radical.

%By \cite{ArdakovBrown} and \cite{Lazard} we know that for any compact
%$p$-adic analytic group $G$, the Iwasawa algebra $\Omega_G$ is a
%complete noetherian semilocal algebra. Combining this fact
%with Corollary \ref{xxsec0.2} we get

\begin{corollary}\label{xxsec1.2}
For any compact $p$-adic analytic group $G$, the Iwasawa algebra
$\Omega_G$ is the dual algebra of an artinian coalgebra $C$.
\end{corollary}

Let $C$ be the dual artinian coalgebra of the Iwasawa algebra
$\Omega_G$. $_{\Omega_G}M$ denotes the category of left
$\Omega_G$-modules and $\Rat(_{\Omega_G}\mathcal{M})$ denotes
the subcategory of $_{\Omega_G}\mathcal{M}$ consisting of all
rational left $\Omega_G$-modules. It is well known that the abelian
category $\mathcal{M}^C$ of right $C$-comodules is equivalent to the
abelian category $\Rat(_{\Omega_G}\mathcal{M})$. Furthermore,
$\Rat(_{\Omega_G}\mathcal{M})$ is closed under submodules,
quotients and arbitrary direct sums. In the sense of \cite{Stenstrom},
it is an hereditary pretorsion class in $_{\Omega_G}\mathcal{M}$.

For the the dual artinian coalgebra $C$ of $\Omega_G$, $C^C$ is an injective
object in $\mathcal{M}^C$, or equivalently, $_{\Omega_G}C$ is an injective
object in ${\rm Rat}(_{\Omega_G}\mathcal{M})$. In general, $_{\Omega_G}C$ is 
not injective in $_{\Omega_G}\mathcal{M}$. But, we have the following result
(see \cite[Section 9.4]{BrzezinskiWisbauer} and \cite[Theorem 3.2]{CuadraNastasescuOystaeyen}). 

\begin{e-proposition}\label{xxsec1.3}
The following statements are equivalent:
\begin{enumerate}
\item[{\rm (a)}] $_{\Omega_G}C$ is injective in $_{\Omega_G}\mathcal{M}$;

\item[{\rm (b)}] $_{\Omega_G}C$ is an injective cogenerator of $_{\Omega_G}\mathcal{M}$;

\item[{\rm (c)}] $C_{\Omega_G}$ is artinian;

\item[{\rm (d)}] $\Omega_G$ is right noetherian;

\item[{\rm (e)}] The injective hull of a rational left $\Omega_G$-module
is rational.
\end{enumerate}
\end{e-proposition}

Throughout this section $C$ is always a left and right artinian coalgebra.
Let $A$ and $B$ be two algebras with its opposite algebras $A^\circ$ and $B^\circ$. 
Let $_AE_B$ be an $(A, B)$-bimodule.
We say that $E$ induces a {\it Morita duality} between $A$ and $B$
if
\begin{enumerate}
\item $_AE$ and $E_B$ are injective cogenerators in the categories
of left $A$-modules and right $B$-modules, respectively; \item the
canonical ring homomorphisms $A\to \End E_B$ and $B^\circ\to \End
{_AE}$ are isomorphisms.
\end{enumerate}
In this case we say that $A$ is {\it left Morita} and $B$ is {\it
right Morita}, and that $A$ is {\it Morita dual to} $B$ (or $A$ and
$B$ are in {\it Morita duality}). If $A=B$, then $A$ is Morita
self-dual, or has a Morita self-duality. We refer to
\cite{AndersonFuller} and \cite{Xue} for some basic properties of a Morita
duality.

For the dual coalgebra $C$ of the Iwasawa algebra $\Omega_G$, let us
consider the $\Omega_G$-bimodule $_{\Omega_G}C_{\Omega_G}$. We easily
observe that $\End(_{\Omega_G}C)\cong \Omega_G^{\circ}$ and 
$\End(C_{\Omega_G})\cong \Omega_G$. Moreover, if $C$ is both
left and right artinian, then $_{\Omega_G}C$ and $C_{\Omega_G}$ are
both injective cogenerators by Proposition \ref{xxsec1.3}. It
follows from \cite[Theorem 24.1]{AndersonFuller} that the
$\Omega_G$-bimodule $_{\Omega_G}C_{\Omega_G}$ defines a Morita
duality. Let $_{\Omega_G}\mathcal{R}$
(resp. $\mathcal{R}_{\Omega_G}$) be the subcategory of
$_{\Omega_G}\mathcal{M}$ (resp. $\mathcal{M}_{\Omega_G}$) consisting
of $_{\Omega_G}C_{\Omega_G}$-reflexive modules. Thus we obtain a
duality
$$
\xymatrix@C=3cm{
_{\Omega_G}\mathcal{R}\ar@<1ex>[r]^{\Hom_{\Omega_G}(- \hspace{2pt},
\hspace{2pt}C)} & \mathcal{R}_{\Omega_G}
\ar@<1ex>[l]^{\Hom_{\Omega_G}(-\hspace{2pt}, \hspace{2pt}C)} }.
\eqno{(\clubsuit)}
$$
Let $_{\Omega_G}\mathcal{M}_{fg}$ be the category of all finitely
generated left $\Omega_G$-modules and $^C\mathcal{M}_{qf}$ be the
category of all quasi-finite left $C$-comodules. Since the Iwasawa
algebra $\Omega_G$ is noetherian and semiperfect, we get $\gld
C=\gld \Omega_G$ by \cite[Proposition
3.4]{CuadraNastasescuOystaeyen}. Applying this fact and the duality
$(\clubsuit)$ yields

\begin{e-proposition}\label{xxsec1.4}
Let $C$ be the dual artinian coalgebra of $\Omega_G$. Then there
exists a duality between abelian categories
$$
\xymatrix@C=3cm{
_{\Omega_G}\mathcal{M}_{fg}\ar@<1ex>[r]^{\Hom_{\Omega_G}(-
\hspace{2pt}, \hspace{2pt}C)} & ^C\mathcal{M}_{qf}
\ar@<1ex>[l]^{(\hspace{5pt})^*}}.
$$
\end{e-proposition}

For any left $\Omega_G$-module $_{\Omega_G}M$, there is a right
$\Omega_G$-module morphism:
$$
\begin{aligned}
\eta_M: \Hom_{\Omega_G}(M, C) & \longrightarrow M^* \\
 f & \longmapsto
\varepsilon\circ f.
\end{aligned}
$$
If $M$ is a finitely generated $\Omega_G$-module, then the right
$\Omega_G$-module $\Hom_{\Omega_G}(M, C)$ is actually a rational
$\Omega_G$-module. Moreover, if we look on $\Hom_{\Omega_G}(M, C)$
as a left $C$-comodule, then it is a quasi-finite comodule. Hence the
image of $\eta_M$ is contained in $\Rat(M^*)$. Therefore we
have a natural transformation:
$$
\eta: \Hom_{\Omega_G}(-, C)\longrightarrow \Rat\circ
(\hspace{5pt})^*, \eqno(\spadesuit)
$$
of functors from $_{\Omega_G}\mathcal{M}_{fg}$ to
$^C\mathcal{M}_{qf}$.

\begin{lemma}\label{xxsec1.5}
The above transformation $\eta$ is a natural isomorphism.
\end{lemma}

If $_{\Omega_G}M$ is a finitely generated free module, then
$\eta_M: \Hom_{\Omega_G}(M, C)\longrightarrow \Rat\circ (M)^*$
is obviously an isomorphism. For a general finitely generated module $M$,
$M$ is finitely presented since $\Omega_G$ is noetherian
$$
\xymatrix{
\underset{{\rm finite}}{\bigoplus} \ar[r] &  \underset{{\rm finite}}{\bigoplus} \ar[r] & M \ar[r] & 0}.
$$
The statement follows from the following commutative diagram:
$$
\xymatrix{
0 \ar[r] & \Hom_{\Omega_G}(M, C) \ar[r] \ar[d]^{\eta_M} & \Hom_{\Omega_G}(\underset{\rm finite}{\bigoplus}\Omega_G, C) \ar[r] \ar[d]^\eta & \Hom_{\Omega_G}(\underset{\rm finite}{\bigoplus}\Omega_G, C)\ar[d]^\eta\\
0 \ar[r] & \Rat\circ(M)^* \ar[r] & \Rat\circ\left(\underset{\rm finite}{\bigoplus}\Omega_G\right)^* \ar[r] & \Rat\circ \left(\underset{\rm finite}{\bigoplus}\Omega_G\right)^*.}
$$

Let $A$ be a noetherian algebra with Jacobson radical $J$ such that $A_0=A/J$ is finite dimensional. Let $_AM$ be an $A$-module. An element $m\in M$ is called a \textit{torsion element} if $J^n m=0$ for all $n\gg 0$. Let us set $\Gamma(M)=\left\{\hspace{2pt} m\in M \hspace{2pt}| \hspace{2pt} m \hspace{2pt}{\rm is}\hspace{2pt}  {\rm a}\hspace{2pt} {\rm torsion}\hspace{2pt} {\rm element}\hspace{2pt} \right\}$. Then $\Gamma(M)$ is a submodule of $M$. In fact, we obtain an additive functor
$$
\Gamma: _A\mathcal{M}\longrightarrow _A\mathcal{M},
$$
by sending an $A$-module to its maximal torsion submodule.
Clearly, $\Gamma$ is a left exact functor. The functor $\Gamma$ has another representation $\Gamma(M)=\underrightarrow{{\rm lim}}\Hom_A(A/J^n, M)$. We
use $\Gamma^\circ$ to denote the torsion functor on the category of
right $A$-modules.

Now let $C$ be the dual coalgebra of Iwasawa algebra $\Omega_G$, and let
$$
\Rat:\hspace{2pt} _{\Omega_G}\mathcal{M}\longrightarrow \hspace{1pt}_{\Omega_G}\mathcal{M}
$$
be the rational functor. Since $C$ is artinian, every finite
dimensional $\Omega_G$-module is rational by \cite[Proposition 3.1.1 and Remarks 3.1.2]{HeynemanRadford}. Therefore, for
a left $\Omega_G$-module, it follows that $\Rat(M)$ is the
sum of all the finite dimensional submodules of $M$. It should be
remarked that the Jacobson radical $J=C_0^\bot$
and $C_0$ is finite dimensional. Thus $\Gamma(M)$ is
also the sum of all the finite dimensional submodules of $M$.
Hence $\Gamma(M)\cong \Rat(M)$. This gives that the functor
$\Gamma$ is naturally isomorphic to the rational functor $\Rat$.
In what follow, we identify the right derived functor $R\Gamma$
with $R{\rm Rat}$.

Let $C$ be the dual coalgebra of the Iwasawa algebra $\Omega_G$.
Then ${\rm Soc}(C)$ is finite dimensional. This implies that
there are only finitely many non-isomorphic simple right (or left)
$C$-comodules. If $^CM$ is quasi-finite, then ${\rm Soc}(M)$
is finite dimensional. Thus $^CM$ is finitely cogenerated. This means that $^C\mathcal{M}_{qf}$ is a
thick subcategory of $^C\mathcal{M}$. So $\mathcal{D}^+
_{qf}(^C\mathcal{M})$, the derived category of bounded below complexes
of left $C$-comodule with quasi-finite cohomology comodules, 
is a full triangulated subcategory of
$\mathcal{D}^+(^C\mathcal{M})$.
Also, since $\Omega_G$ is noetherian, $\mathcal{D}^-_{
fg}(_{\Omega_g}\mathcal{M})$, the derived category of bounded above
complexes of left $\Omega_G$-modules with finitely generated
cohomology modules, is a full triangulated subcategory of
$\mathcal{D}^-(_{\Omega_G}\mathcal{M})$. The duality in Proposition
\ref{xxsec1.4} gives rise to a duality of derived categories.

\begin{theorem}\label{xxsec1.6}
Let $C$ be the artinian coalgebra of the Iwasawa algebra
$\Omega_G$. We have dualities of triangulated categories:
$$
\xymatrix@C=2cm{
\mathcal{D}^-_{fg}(_{\Omega_G}\mathcal{M})\ar@<1ex>[r]^{R\Gamma^\circ\circ (\hspace{5pt})^*} & \mathcal{D}^+_{qf}(^C\mathcal{M})
\ar@<1ex>[l]^{(\hspace{5pt})^*}},\hspace{8pt} \xymatrix@C=2cm{
\mathcal{D}^b_{fg}(_{\Omega_G}\mathcal{M})\ar@<1ex>[r]^{R\Gamma^\circ\circ (\hspace{5pt})^*} & \mathcal{D}^b_{qf}(^C\mathcal{M})
\ar@<1ex>[l]^{(\hspace{5pt})^*}}
$$
\end{theorem}

By the fact that $C$ is artinian, we know that $\mathcal{D}^+_{qf}(^C\mathcal{M})$
is equivalent to the derived category $\mathcal{D}^+(^C\mathcal{M}_{qf})$
of complexes of quasi-finite comodules. Since $\Omega_G$ is notherian, $\mathcal{D}^-_{fg}(_{\Omega_G}\mathcal{M})$ is equivalent to $\mathcal{D}^-(_{\Omega_G}\mathcal{M}_{fg})$. In view of Proposition \ref{xxsec1.4}
we obtain the following duality
$$
\xymatrix@C=3cm{
\mathcal{D}^-(_{\Omega_G}\mathcal{M}_{fg})\ar@<1ex>[r]^{\Hom_{\Omega_G}(-,\hspace{5pt} C)} & \mathcal{D}^+(^C\mathcal{M}_{qf})
\ar@<1ex>[l]^{(\hspace{5pt})^*}}.
$$
Applying Lemma \ref{xxsec1.5} we easily verify that the composition
\[
\mathcal{D}^-_{fg}(_{\Omega_G}\mathcal{M})\xrightarrow[]{\cong}  \mathcal{D}^-(_{\Omega_G}\mathcal{M}_{fg}) \xrightarrow[]{\Hom_{\Omega_G}(-,C)} \mathcal{D}^+(^C\mathcal{M}_{qf}) \xrightarrow[]{\cong}  \mathcal{D}^+_{qf}(^C\mathcal{M})
\]
is isomorphic to the functor $R\Gamma^\circ \circ (\hspace{5pt})^*$. Furthermore, one can observe that $R\Gamma^\circ \circ (\hspace{5pt})^*$ sends bounded complexes to
bounded complexes.

\begin{corollary}\label{xxsec1.7}
Let $C$ be the dual artinian coalgebra of the Iwasawa algebra $\Omega_G$.
Then we have a duality of triangulated categories:
$$
\xymatrix@C=3cm{
\mathcal{D}^b_{fd}(_{\Omega_G}\mathcal{M})\ar@<1ex>[r]^{R\Gamma\circ( \hspace{5pt} )^*} & \mathcal{D}^b_{fd}(^C\mathcal{M})
\ar@<1ex>[l]^{(\hspace{5pt})^*}}.
$$
\end{corollary}

It is enough to show that $R\Gamma\circ (M)^*\in \mathcal{D}^b_{fd}(^C\mathcal{M})$
for any finite dimensional left $\Omega_G$-module $M$. Since $\Omega_G$
is noetherian and complete with respect to the $J$-adic filtration,
the Jacobson radical $J$ of $\Omega_G$ satisfies Artin-Rees condition.
Applying \cite[Theorem 3.2]{CuadraNastasescuOystaeyen} yields that
the injective envelop of a $J$-torsion module is still $J$-torsion.
Now $M^*$ is a $J$-torsion module. Thus we have an injective resolution of $M^*$
with each component being $J$-torsion. Therefore $R\Gamma(M^*)$
is quasi-isomorphic to $M^*$, that is, $R\Gamma(M^*)\in \mathcal{D}^b_{fd}(^C\mathcal{M})$.

\section{Conslusion}
\label{xxsec2}

It is well-known that there is an informal idea in the theory of derived categories of coalgebras: `the derived
category associated to a coalgebra keeps the homological relevant information, so
if two coalgebras have equivalent derived categories, they should have isomorphic
cohomology'. In view of the beautiful structure and properties of the Iwasawa algebra $\Omega_G$,
$C$ as the dual artinian coalgebra of $\Omega_G$ should have elegant properties which 
build stable basis for the further study of derived categories of bounded complexes of 
$C$-comodules. Let $D^b_{qf}(^C\mathcal{M})$ (respectively, $D^b_{qf}(\mathcal{M}^C)$) be the derived category of bounded complexes
of left (respectively, right) $C$-comodules with quasi-finite cohomology comodules.
We wonder whether it is possible for us to establish a duality of the derived categories of left $C$-comodules and of right $C$-comodules by the results obtained in Section \ref{xxsec1}. More precisely speaking, we have the following question:

\begin{question}
Let $C$ be the dual artinian coalgebra of the Iwasawa algebra $\Omega_G$.
Are there dualities of triangulated categories
$$
\xymatrix{
\mathcal{D}^b_{qf}(^C\mathcal{M})\ar@<1ex>[r]^F &
\mathcal{D}^b_{qf}(\mathcal{M}^C)
\ar@<1ex>[l]^G } \hspace{7pt}?
$$
Here $F=R\Gamma\circ(\hspace{5pt})^*$ and $G=R\Gamma^\circ \circ (\hspace{5pt})^*$.
\end{question}

\bigskip

\noindent {\bf Acknowledgements} I have accumulated quite a debt of
gratitude in writing this paper. significantly help us improve the final presentation of this
%article.

\end{document}